\documentclass[12pt,oneside,reqno]{amsart}
\usepackage[all]{xy}
\usepackage{amsfonts,amsmath,oldgerm,amssymb,amscd}
\UseComputerModernTips
\numberwithin{equation}{section}
\usepackage[breaklinks]{hyperref}
\allowdisplaybreaks

\usepackage{enumerate}
\usepackage{caption} 
\newtheorem{theorem}{Theorem}[section]

\newtheorem{lemma}[subsection]{{\bf Lemma}}

\begin{document}

	\title[Balancing and Lucas-balancing numbers as difference of two repdigits]{Balancing and Lucas-balancing numbers as difference of two repdigits} 

\author[M. Mohapatra]{M. Mohapatra}
\address{Monalisa Mohapatra, Department of Mathematics, National Institute of Technology Rourkela, Odisha-769 008, India}
\email{mmahapatra0212@gmail.com}

\author[P. K. Bhoi]{P. K.  Bhoi}
\address{Pritam Kumar Bhoi, Department of Mathematics, National Institute of Technology Rourkela, Odisha-769 008, India}
\email{pritam.bhoi@gmail.com}

\author[G. K. Panda]{G. K. Panda}
\address{Gopal Krishna Panda, Department of Mathematics, National Institute of Technology Rourkela, Odisha-769 008, India}
\email{gkpanda\_nit@rediffmail.com}

\thanks{2020 Mathematics Subject Classification: Primary 11B39, Secondary 11J86, 11D61. \\
	Keywords: balancing sequence, Lucas-balancing sequence, linear forms in logarithms, Baker-Davenport reduction method}

\begin{abstract}
	Positive integers with all digits equal are called repdigits. In this paper, we find all balancing and Lucas-balancing numbers, which can be expressed as the difference of two repdigits. The method of proof involves the application of Baker's theory for linear forms in logarithms of algebraic numbers and the Baker-Davenport reduction procedure.
\end{abstract}

\maketitle
\pagenumbering{arabic}
\pagestyle{headings}

	\section{Introduction}

The balancing number sequence $(B_n)_{n\geq0}$ and the Lucas-balancing sequence $(C_n)_{n\geq0}$ are defined by the binary recurrences
\begin{equation}\label{bn}
	B_{n+1}=6B_n-B_{n-1},\quad B_0 = 0,\quad B_1 = 1
\end{equation}
and 
\begin{equation}\label{cn}
	C_{n+1} = 6C_{n}-C_{n-1},\quad C_0=1,\quad C_1=3.
\end{equation}
The Binet’s formula for the sequences are given by
\begin{center}
	$B_n=\frac{\alpha^{n}-\beta^{n}}{4\sqrt{2}}$
	and 
	$C_n=\frac{\alpha^{n}-\beta^{n}}{2};$ for $n\geq1$,
\end{center}
where $(\alpha,\beta)=(3+2\sqrt{2},3-2\sqrt{2})$ is the pair of roots of the characteristic polynomial $x^2-6x-1$. This implies easily that the inequalities
\begin{center}
	$\alpha^{n-1} \leq B_n <\alpha^n$ and $\alpha^{n} < 2C_n <\alpha^{n+1}$,
\end{center} 
hold for all $n\geq1$.\\

A repdigit is a positive integer in which all digits are the same. The study of repdigits within second-order linear recurrence sequences has garnered interest among mathematicians. In \cite{rp2018}, all repdigit instances of Balancing and Lucas-balancing numbers were identified. S. G. Rayaguru and G. K. Panda \cite{rp2021} explored Balancing and Lucas-balancing numbers that can be expressed as the sums of two repdigits. Mohapatra et al. \cite{mbp2024} examined the existence of repdigits as the difference between two Balancing or Lucas-balancing numbers. Erduwan et al. \cite{ekl2021} identified all Fibonacci and Lucas numbers that can be represented as the difference of two repdigits. Edjeou and Faye \cite{ef2023} accomplished a similar task for Pell and Pell-Lucas numbers. More recently, M.G. Duman \cite{d2023} found all Padovan numbers that can be expressed as the difference of two repdigits.

In this paper, we investigate a related question: which Balancing and Lucas-balancing numbers can be written as the difference of two repdigits? To address this, we consider the following two equations:
\begin{equation}\label{1}
	B_k = d_1\biggr(\frac{10^n -1}{9}\biggr)- d_2\biggr(\frac{10^m -1}{9}\biggr)
\end{equation}
and
\begin{equation}\label{2}
	C_k = d_1\biggr(\frac{10^n -1}{9}\biggr)- d_2\biggr(\frac{10^m -1}{9}\biggr),
\end{equation}
where $(k, m, n)$ are positive integers with $n\geq2$.
\section{Auxiliary Results}
To solve the Diophantine equations involving repdigits and terms of binary recurrence sequences, many authors have used Baker’s theory to reduce lower bounds concerning linear forms in logarithms of algebraic numbers. These lower bounds play an important role in solving such Diophantine equations. We start by recalling some basic definitions and results from algebraic number theory.

Let $\lambda$ be an algebraic number with minimal primitive polynomial 
\[f(X) = a_0(X-\lambda^{(1)})\cdot\cdot\cdot(X-\lambda^{(k)}) \in \mathbb{Z}[X],\]
where $a_0 > 0$ and $\lambda^{(i)}$'s are conjugates of $\lambda$. Then
\begin{center}
	$h(\lambda) = \dfrac{1}{k}\biggr(\log a_0 + \sum\limits_{j=1}^{k}$max$\{0,\log|\lambda^{(j)}|\}\biggr)$
\end{center}
is called the $logarithmic$ $height$ of $\lambda$. If $\lambda = \dfrac{a}{b}$ is a rational number with $gcd(a,b)=1$ and $b>1$, then  $h(\lambda)=\log$(max$\{|a|,b\})$. We give some properties of the logarithmic height whose proofs can be found in \cite{bms2006}.

Let $\gamma$ and $\eta$ be two algebraic numbers. then 
\begin{equation*}
	\begin{split}
		& (i)\hspace{0.2cm} h(\gamma\pm\eta) \leq h(\gamma)+h(\gamma)+\log 2,\\
		&(ii)\hspace{0.2cm} h(\gamma\eta^{\pm 1})\leq h(\gamma)+h(\eta),\\
		&(iii)\hspace{0.2cm} h(\gamma ^k)=|k|h(\gamma).
	\end{split}
\end{equation*}
Now, we give a theorem which is  deduced from  Corollary  2.3  of  Matveev \cite{m2000} and provides a large upper bound for the subscript $k$ in the Equation \eqref{1} and Equation \eqref{2} (also see Theorem 9.4 in \cite{bms2006}).

\begin{theorem}\label{thm1}\cite{m2000}. Let $\mathbb{L}$ be an algebraic number field of degree $d_{\mathbb{L}}$. Let $\gamma_1,\ldots,\gamma_l \in \mathbb{L}$ be positive real numbers and $b_1, \ldots, b_l$ be nonzero integers. If $\Gamma = \prod\limits_{i=1}^{l} \gamma_{i}^{b_i} -1$ is not zero, then 
	\begin{equation*}
		\log|\Gamma| > -1.4\cdot 30^{l+3}\cdot l^4.5 \cdot d_{\mathbb{L}}^2 (1+\log d_{\mathbb{L}}^2)(1+\log(D))A_1 A_2 \cdot\cdot\cdot A_l,
	\end{equation*}
	where $D\geq$ max$\{|b_1|,\ldots,|b_l|\}$ and  $A_1, \cdots, A_l$ are positive integers such that $A_j \geq h'(\gamma_j)$ = $max\{d_{\mathbb{L}}h(\gamma_j),|\log \gamma_j|, 0.16\},$ for $j= 1,\ldots,l$.
\end{theorem}

The following result of Baker-Davenport due to Dujella and Peth\H{o} \cite{dp1998} is another tool in our proofs. It will be used to reduce the upper bounds of the variables on Equation \eqref{1} and Equation \eqref{2}.
\newline

\begin{lemma}\label{lem1}\cite{dp1998}. Let $M$ be a positive integer and $\frac{p}{q}$ be a convergent of the continued fraction of the irrational number $\tau$ such that $q > 6M$. Let $A, B, \mu$ be some real numbers with $A > 0$ and $B > 1$. Let $\epsilon := \left\lVert \mu q\right\rVert - M\left\lVert \tau q\right\rVert$, where $\left\lVert .\right\rVert$ denotes the distance from the nearest integer. If $\epsilon > 0$, then there exists no solution to the inequality
	\begin{equation}
		0 <|u\tau - v +\mu|<AB^{-w},
	\end{equation}
	in positive integers $u,v,w$ with $u\leq M$ and $w\geq \dfrac{\log(Aq/\epsilon)}{\log B}$.
\end{lemma}

\begin{lemma}\label{lem2}\cite{wb1989}. 
	Let $a,x\in \mathbb R$. If $0<a<1$ and $|x|<a$, then 
	\begin{equation*}
		|\log(1+x)|<\dfrac{-|log(1-a)|}{a}.|x|
	\end{equation*} and
	\begin{equation*}
		|x|<\dfrac{a}{1-e^{-a}}.|e^x - 1|.
	\end{equation*}
\end{lemma}
\begin{lemma}\cite{rp2020}\label{lem3}.	The only balancing number which is concatenation of two repdigits is 35.
\end{lemma}
\begin{lemma}\cite{rb2023}\label{lem4}.	The only balancing numbers which are concatenation of three repdigits are 204 and 1189.
\end{lemma}
\begin{lemma}\cite{rps2021}\label{lem5}.	The only Lucas-balancing numbers which are concatenation of two repdigits are 17 and 577.
\end{lemma}
\begin{lemma}\cite{rb2023}\label{lem6}.	The only Lucas-balancing number which is concatenation of three repdigits is 3363.
\end{lemma}
\section{Balancing numbers as difference of two repdigits}
\begin{theorem}\label{thm2} If $B_k$ is expressible as difference of two repdigits, then
	$B_k \in \{6, 35\}$, where
	$B_2 = 6 = 11-5$ and
	$B_3 = 35 = 44-9$.    
	\begin{proof} Assume that Equation \eqref{1} holds. Let $1 \leq k \leq 25$ and $n\geq 2$. Then, by using $Mathematica$,
		we obtain only the solutions listed in Theorem \ref{thm2}.
		From now on, assume that $k > 25$. If $n=m$, then it follows that $d_1 > d_2$, which means that $B_k$ is a repdigit. However, the largest repdigit in $B_k$ is 6 \cite{rp2018}. Thus, we get a contradiction since $k > 25$. If $n-m=1$ and $d_1\geq d_2$, we encounter balancing numbers that are concatenations of two repdigits. This scenario is impossible, as established in Lemma \ref{lem3}. If $n-m=1$ and $d_1< d_2$, then we have balancing numbers that are concatenation of three repdigits. This is impossible by Lemma \ref{lem4}. Then we may suppose that $n - m \geq 2$. The inequality
		\begin{equation*}
			\dfrac{\alpha^n}{20}<\dfrac{10^{n-1}}{2}<10^{n-1}-10^{m-1}<\dfrac{d_1(10^n-1)}{9}-\dfrac{d_2(10^m-1)}{9}=B_k<\alpha^k,
		\end{equation*}
		implies that $n<k+2$. On the other hand, we can elegantly rewrite Equation \eqref{1} as
		\begin{equation*}
			\dfrac{\alpha^k-\beta^k}{4\sqrt{2}}= \dfrac{d_1(10^n -1)}{9}- \dfrac{d_2(10^m -1)}{9}.
		\end{equation*}
		This transformation leads us to the refined expression:
		\begin{equation}\label{3}
			\dfrac{9\alpha^k}{4\sqrt{2}}- d_1 10^n =  \dfrac{9\beta^k}{4\sqrt{2}}- d_2 10^m -(d_1-d_2).
		\end{equation}
		Taking absolute value of both sides of the equality \eqref{3}, we obtain
		\begin{equation}\label{4}
			\biggr|\dfrac{9\alpha^k}{4\sqrt{2}}-d_1 10^n\biggr| \leq \dfrac{9|\beta|^k}{4\sqrt{2}}+d_2 10^m + |d_1-d_2|.
		\end{equation}
		Dividing both sides of \eqref{4} by $d_1 10^n$, we arrive at
		\begin{equation*}
			\begin{split}
				\biggr|\dfrac{9\cdot10^{-n}\cdot\alpha^k}{d_1\cdot 4\sqrt{2}}-1\biggr| & \leq \dfrac{9|\beta|^k}{d_1 10^n\cdot 4\sqrt{2}} + \dfrac{d_2 10^m}{d_1 10^n}+\dfrac{|d_1-d_2|}{d_1 10^n} \\ 
				& \leq \dfrac{9|\beta|^k}{10^{n-m+1}\cdot4\sqrt{2}} + \dfrac{9}{10^{n-m}}+\dfrac{8}{10^{n-m+1}} . \\ 
			\end{split}
		\end{equation*}
		From this, we conclude
		\begin{equation} \label{5}
			\biggr|\dfrac{9\cdot10^{-n}\cdot\alpha^k}{d_1\cdot4\sqrt{2}}-1\biggr| < \dfrac{9.81}{10^{n-m}}.
		\end{equation}
		We now invoke Theorem \ref{thm1} with $(\gamma_1, \gamma_2, \gamma_3)$ = $(\alpha,10, 9/(d_1 \cdot 4\sqrt{2}))$ and $(b_1, b_2, b_3)$ = $(k, -n, 1)$. Note that $\gamma_1$, $\gamma_2$, and $\gamma_3$ are positive real numbers and elements of the field $\mathbb K = \mathbb Q(\sqrt{2})$. Therefore, the degree of the field
		$\mathbb K$ is equal to $d_\mathbb{L} = 2$. Let us define
		\begin{equation*}
			\Gamma_1 = \dfrac{9\cdot10^{-n}\cdot\alpha^k}{d_1\cdot 4\sqrt{2}}-1 .
		\end{equation*}
		Assume that $\Gamma_1=0$. Then $\alpha^k = \dfrac{10^n d_1 \cdot 4\sqrt{2}}{9}$. Conjugating in $\mathbb Q(\sqrt{2})$ gives us $\beta^k = -10^n d_1 \cdot 4\sqrt{2}/9$ and so $C_k = \alpha^k + \beta^k = 0$, which presents a contradiction. Hence, we affirm that $\Gamma_1 \neq 0$.

		Utilizing the properties of absolute logarithmic height, we determine the logarithmic heights of $\gamma_1$, $\gamma_2$ and $\gamma_3$ as $h(\gamma_1)= \frac{\log\alpha}{2}$, $h(\gamma_2)= \log 10$ and $h(\gamma_3)\leq h(4d_1\sqrt{2})+h(9) < 6.2$.
		
		We can take $A_1=\log \alpha$, $A_2= 2\log 10$, $A_3= 12.4$. Since $n < k+2$, we take $D$ = max$\{n, k, 1\} = k+2$. 
		With \eqref{5} in mind and by invoking Theorem \ref{thm1}, we arrive at
		\begin{equation*}
			\log\biggr(\dfrac{9.81}{10^{n-m}}\biggr)>\log|\Gamma_1| > -1.4\cdot30^6\cdot3^{4.5}\cdot2^2 (1+\log 2)(1+\log n)(\log\alpha)(2\log 10)(12.4).
		\end{equation*}
		By a simple calculation, the above inequality leads to 
		\begin{equation}\label{6}
			\begin{split}
				(m-n)\log(\alpha) &< \log (9.81) + 9.8 \cdot 10^{13} (1+\log (k+2))\\ & < 9.9 \cdot 10^{13} (1+\log (k+2)).\\
			\end{split}
		\end{equation}
		Rearranging Equation \eqref{1} as 
		\begin{equation} \label{7}
			\dfrac{\alpha^k}{4\sqrt{2}}- \dfrac{d_1 10^n - d_2 10^m}{9} = \dfrac{\beta^k}{4\sqrt{2}}- \dfrac{(d_1 - d_2)}{9}
		\end{equation}
		and taking absolute value of both sides of Equation \eqref{7}, we get
		\begin{equation}\label{8}
			\biggr|\dfrac{\alpha^k}{4\sqrt{2}}- \dfrac{d_1 10^n - d_2 10^m}{9}\biggr| = \dfrac{|\beta^k|}{4\sqrt{2}}- \dfrac{|d_1 - d_2|}{9}.
		\end{equation}
		Dividing both sides of the above inequality by $\frac{\alpha^k}{4\sqrt{2}}$, we obtain
		\begin{equation}\label{9}
			\biggr|1- \dfrac{(d_1-d_210^{m-n})\cdot10^n\cdot\alpha^{-k}}{9\cdot 4\sqrt{2}}\biggr| \leq \frac{1}{\alpha^{2k}}+\dfrac{8\cdot 4\sqrt{2}}{9\alpha^k}< \frac{6}{\alpha^k}.
		\end{equation}
		Now, we can again apply Theorem \ref{thm1} to the above inequality with
		\begin{center}
			$(\gamma_1, \gamma_2, \gamma_3)$ = $\biggr(\alpha, 10, \dfrac{(d_1-d_210^{m-n})}{9\cdot4\sqrt{2}}\biggr)$ and $(b_1, b_2, b_3)$ = $(k,-n, 1)$.
		\end{center}
		Note that $\gamma_1$, $\gamma_2$, and $\gamma_3$ are positive real numbers and elements of the field $\mathbb K = \mathbb Q(\sqrt{2})$. Therefore, the degree of the field $\mathbb K$ is equal to $d_\mathbb{L} = 2$.
		Let
		\begin{equation*}
			\Gamma_2 = 1- \dfrac{(d_1-d_210^{m-n})\cdot10^n\cdot\alpha^{-k}}{9\cdot 4\sqrt{2}}.
		\end{equation*}
		If $\Gamma_2=0$, then $\alpha^{2k}\in \mathbb Q$, which is false for $k>0$. By using the properties of the absolute $logarithmic$ $height$, we get
		\begin{equation*}
			h(\gamma_1)=h(\alpha)=\dfrac{\log\alpha}{2},\hspace{0.2cm}
			h(\gamma_2)= \log 10, 
		\end{equation*}
		and 
		\begin{equation*}
			\begin{split}
				h(\gamma_3) &=  h\biggr(\dfrac{(d_1-d_210^{m-n})}{9\cdot4\sqrt{2}}\biggr)\\ &\leq h(9)+ h(4\sqrt{2})+h(d_1)+h(d_2)+(n-m)\log10 +\log 2 \\ & < 9.02+(n-m) \log 10.
			\end{split}
		\end{equation*}
		So we can take $A_1=\log \alpha$, $A_2= 2\log 10$, $A_3= 18.04+2(n-m)\log 10$. Since $n < k+2$, we take $D$ = max$\{n, k, 1\} = k+2$. 
		Thus, taking into account \eqref{9} and using Theorem \ref{thm1}, we obtain
		\begin{equation*}
			6\cdot\alpha^k>|\Gamma_2|> e^{(C\cdot(1+\log 2)(1+\log(k+2))\cdot\log \alpha\cdot2\log 10\cdot(18.04+2(n-m)\log 10))},
		\end{equation*}
		where $C= -1.4\cdot30^6\cdot3^{4.5}\cdot 2^2$. By a simple computation, it follows that 
		\begin{equation}\label{10}
			k\log\alpha-\log 6 < 7.9\cdot10^{12}\cdot(1+\log(k+2))(18.04+2(n-m)\log 10).
		\end{equation}
		Using \eqref{6} and \eqref{10}, a computer search with $Mathematica$ gives us that $k < 1.7\cdot10^{28}$.\\
		
		Now, let us reduce the upper bound on $n$ by using the Baker–Davenport algorithm as outlined in Lemma \ref{lem1}. Let
		\begin{equation*}
			\Lambda_1=k\log\alpha - n\log10 +\log\biggr(\dfrac{9}{d_1\cdot 4\sqrt{2}}\biggr).
		\end{equation*}
		
		From \eqref{5}, we have
		\begin{equation*}
			|x| = |e^{\Lambda_1}-1|< \dfrac{9.81}{10^n-m}<\frac{1}{10}
		\end{equation*}
		for $n-m \geq 2$. Choosing $a = 0.1$, we get the inequality
		\begin{equation*}
			|\Lambda_1|=|\log(x+1)|< \dfrac{\log(10/9)}{1/10} \cdot \frac{9.81}{10^{n-m}}<(10.34)\cdot 10^{m-n}
		\end{equation*}
		by Lemma \ref{lem2}. Thus, it follows that 
		\begin{equation*}
			0<\biggr|k\log\alpha - n\log10 +\log\biggr(\dfrac{9}{d_1\cdot 4\sqrt{2}}\biggr)\biggr|< (10.34)\cdot 10^{m-n}.
		\end{equation*}
		Dividing this inequality by $\log 10$, we obtain
		\begin{equation}\label{11}
			0<\biggr|k\biggr(\frac{\log\alpha}{\log 10}\biggr) - n +\dfrac{\log\biggr(\dfrac{9}{d_1\cdot 4\sqrt{2}}\biggr)}{\log 10}\biggr|< (4.5) \cdot 10^{m-n}.
		\end{equation}
		We can take $\tau=\frac{\log\alpha}{\log 10} \notin \mathbb Q$ and $M = 1.7\cdot 10^{28}$.\\ Then we found that $q_{59} = 808643106803003389273254071835$, the denominator of the $59^{th}$ convergent of $\tau$ is greater than $6M$.
		Now take $\mu=\dfrac{\log(9/d_1\cdot 4\sqrt{2})}{\log 10}$.
		In this case, considering the fact that $1 \leq d_1 \leq9$, a quick computation with $Mathematica$ gives us that $\epsilon(\mu) := \left\lVert \mu q_{59}\right\rVert - M\left\lVert \tau q_{59}\right\rVert = 0.03855$. Let $A=4.5$, $B= 10$, and $\omega=n-m$ in Lemma \ref{lem1}. Thus, using $Mathematica$, we can say that \eqref{11} has no solution if 
		\begin{equation*}
			\dfrac{\log(Aq_{59}/\epsilon(\mu))}{\log B} < 31.97<n-m.
		\end{equation*}
		So $n-m \leq 31$.
		Substituting this upper bound for $n - m$ in \eqref{10}, we obtain $k < 3.1 \cdot 10^{15} $. Now, let
		\begin{equation*}
			\Lambda_2= n\log10-k\log\alpha +\log\biggr(\dfrac{d_1-d_210^{m-n}}{9\cdot 4\sqrt{2}}\biggr).
		\end{equation*}
		From \eqref{9}, we have 
		\begin{equation*}
			|x| = |e^{\Lambda_2}-1|< \dfrac{6}{\alpha^k}<\frac{1}{10}
		\end{equation*}
		for $k\geq 25$. Choosing $a = 0.1$, we get the inequality
		\begin{equation*}
			|\Lambda_2|=|\log(x+1)|< \dfrac{\log(10/9)}{1/10} \cdot \frac{6}{\alpha^k}<(6.33)\cdot \alpha^{-k}
		\end{equation*}
		by Lemma \ref{lem2}. Thus, it follows that
		\begin{equation*}
			0<\biggr|n\log10-k\log\alpha +\log\biggr(\dfrac{d_1-d_210^{m-n}}{9\cdot 4\sqrt{2}}\biggr)|< (6.33)\cdot \alpha^{-k}.
		\end{equation*}
		Dividing both sides of the above inequality by $\log \alpha$, we obtain
		\begin{equation}\label{12}
			0<\biggr|n\biggr(\frac{\log 10}{\log \alpha}\biggr) - k +\dfrac{\log\biggr(\dfrac{d_1-d_210^{m-n}}{9\cdot 4\sqrt{2}}\biggr)}{\log \alpha}\biggr|< (3.6)\cdot \alpha^{-k}.
		\end{equation}
		Putting $\tau = \frac{\log 10}{\log \alpha} \notin \mathbb Q$ and taking $M= 3.1\cdot10^{15}$, we found that $q_{36} = 73257846218558279$, the denominator of the $36^{th}$ convergent of $\tau$ exceeds $6M$. Now take
		$\mu = \dfrac{\log\biggr(\dfrac{d_1-d_210^{m-n}}{9\cdot 4\sqrt{2}}\biggr)}{\log \alpha}$.
		
		In this case, considering the fact that $1 \leq d_1$, $d_2 \leq 9$ and $2 \leq n - m \leq 31$, a quick computation gives us the inequality $\epsilon(\mu) := \left\lVert \mu q_{36}\right\rVert - M\left\lVert \tau q_{36}\right\rVert = 0.327932$.
		
		Let $A=3.6$, $B= \alpha$, and $\omega=k$ in Lemma \ref{lem1}. Thus, using $Mathematica$, we can say that \eqref{12} has no solution if 
		\begin{equation*}
			\dfrac{\log(Aq_{36}/\epsilon(\mu))}{\log B} < 23.38< k.
		\end{equation*}
		This implies $k \leq 23$, which contradicts our assumption that $k\geq 25$. Thus, the proof is completed.

	\end{proof}
\end{theorem}

\section{Lucas-balancing numbers as difference of two repdigits}

\begin{theorem}\label{thm3} If $C_k$ is expressible as difference of two repdigits, then 
	$C_k \in \{3, 17\}$, where
	$C_2 = 3 = 11-8$ and
	$C_3 = 17 = 22-5$.
	\begin{proof} Assume that \eqref{2} holds. Let $1 \leq k \leq 25$ and $n\geq 2$. Then, by using $Mathematica$,
		we obtain only the solutions listed in Theorem \ref{thm3}.
		From now on, assume that $k \geq 25$. If $n=m$, then it follows that $d_1 > d_2$, which means that $C_k$ is a repdigit. But the largest repdigit in $C_k$ is 99 \cite{rp2018}. Thus, we get a contradiction since $k \geq 25$. If $n-m=1$ and $d_1\geq d_2$, then we have Lucas-balancing numbers that are concatenation of two repdigits. This is impossible by Lemma \ref{lem5}. If $n-m=1$ and $d_1< d_2$, then we have Lucas-balancing numbers that are concatenation of three repdigits. This is impossible by Lemma \ref{lem6}. Then we may suppose that $n - m \geq 2$. The inequality
		\begin{equation*}
			\dfrac{\alpha^n}{20}<\dfrac{10^{n-1}}{2}<10^{n-1}-10^{m-1}<\dfrac{d_1(10^n-1)}{9}-\dfrac{d_2(10^m-1)}{9}=C_k<\alpha^{k+1},
		\end{equation*}
		implies that $n<k+3$. On the other hand, we rewrite \eqref{2} as 
		\begin{equation*}
			\dfrac{\alpha^k-\beta^k}{2}= \dfrac{d_1(10^n -1)}{9}- \dfrac{d_2(10^m -1)}{9}
		\end{equation*}
		to obtain
		\begin{equation}\label{13}
			\dfrac{9\alpha^k}{2}- d_1 10^n =  \dfrac{9\beta^k}{2}- d_2 10^m -(d_1-d_2).
		\end{equation}
		Taking absolute value of both sides of \eqref{13}, we obtain
		\begin{equation}\label{14}
			\biggr|\dfrac{9\alpha^k}{2}-d_1 10^n\biggr| \leq \dfrac{9|\beta|^k}{2}+d_2 10^m + |d_1-d_2|.
		\end{equation}
		Dividing both sides of \eqref{14} by $d_1 10^n$, we obtain
		\begin{equation*}
			\begin{split}
				\biggr|\dfrac{9\cdot10^{-n}\cdot\alpha^k}{2d_1}-1\biggr| & \leq \dfrac{9|\beta|^k}{2d_1 10^n} + \dfrac{d_2 10^m}{d_1 10^n}+\dfrac{|d_1-d_2|}{d_1 10^n} \\ 
				& \leq \dfrac{9|\beta|^k}{2\cdot10^{n-m+1}} + \dfrac{9}{10^{n-m}}+\dfrac{8}{10^{n-m+1}} . \\ 
			\end{split}
		\end{equation*}
		From this, it follows that
		\begin{equation} \label{15}
			\biggr|\dfrac{9\cdot10^{-n}\cdot\alpha^k}{2d_1}-1\biggr| < \dfrac{9.81}{10^{n-m}}.
		\end{equation}
		We now apply Theorem \ref{thm1} with $(\gamma_1, \gamma_2, \gamma_3)$ = $(\alpha, 10, 9/2d_1)$ and $(b_1, b_2, b_3)$ = $(k, -n, 1)$. Note that $\gamma_1$, $\gamma_2$, and $\gamma_3$ are positive real numbers and elements of the field $\mathbb K = \mathbb Q(\sqrt{2})$. Therefore, the degree of the field
		$\mathbb K$ is equal to $D = 2$. Put
		\begin{equation*}
			\Gamma_3 = \dfrac{9\cdot10^{-n}\cdot\alpha^k}{2d_1}-1 .
		\end{equation*}
		Assume that $\Gamma_3=0$. Then $\alpha^{k}= \dfrac{2d_110^n}{9}\in \mathbb Q$, which is false for $k>0$. Therefore, $\Gamma_3 \neq 0$.

		By the properties of the absolute $logarithmic$ $height$, the logarithmic heights of $\gamma_1$, $\gamma_2$ and $\gamma_3$ are calculated as $h(\gamma_1)=\frac{\log\alpha}{2}$, $h(\gamma_2)= \log 10$ and $h(\gamma_3)\leq h(2d_1)+h(9) < 5.09 $.
		
		We can take $A_1=\log \alpha$, $A_2= 2\log 10$, $A_3= 10.18$. Since $n < k+3$, we take $D$ = max$\{n, k, 1\} = k+3$. 
		Thus, taking into account \eqref{15} and using Theorem \ref{thm1}, we obtain
		\begin{equation*}
			\log\biggr(\dfrac{9.81}{10^{n-m}}\biggr)>\log|\Gamma_3| > -1.4\cdot30^6\cdot3^{4.5}\cdot2^2 (1+\log 2)(1+\log n)(\log\alpha)(2\log 10)(10.18).
		\end{equation*}
		By a simple calculation, the above inequality leads to 
		\begin{equation}\label{16}
			\begin{split}
				(m-n)\log(\alpha) &< \log (9.81) + 8.02 \cdot 10^{13} (1+\log (k+3))\\ & < 8.03 \cdot 10^{13} (1+\log (k+3)).\\
			\end{split}
		\end{equation}
		Rearranging \eqref{2} as 
		\begin{equation} \label{17}
			\dfrac{\alpha^k}{2}- \dfrac{d_1 10^n - d_2 10^m}{9} = \dfrac{\beta^k}{2}- \dfrac{(d_1 - d_2)}{9}
		\end{equation}
		and taking absolute value of both sides of \eqref{17}, we get
		\begin{equation}\label{18}
			\biggr|\dfrac{\alpha^k}{2}- \dfrac{d_1 10^n - d_2 10^m}{9}\biggr| = \dfrac{|\beta^k|}{2}- \dfrac{|d_1 - d_2|}{9}.
		\end{equation}
		Dividing both sides of the above inequality by $\frac{\alpha^k}{2}$, we obtain
		\begin{equation}\label{19}
			\biggr|1- \dfrac{(d_1-d_210^{m-n})\cdot10^n\cdot\alpha^{-k}}{18}\biggr| \leq \frac{1}{\alpha^{2k}}+\dfrac{16}{9\alpha^k}< \frac{3}{\alpha^k}.
		\end{equation}
		Now, we can again apply Theorem \ref{thm1} to the above inequality with
		\begin{center}
			$(\gamma_1, \gamma_2, \gamma_3)$ = $\biggr(\alpha, 10, \dfrac{(d_1-d_210^{m-n})}{18}\biggr)$ and $(b_1, b_2, b_3)$ = $(k, -n, 1)$.
		\end{center}
		Note that $\gamma_1$, $\gamma_2$, and $\gamma_3$ are positive real numbers and elements of the field $\mathbb K = \mathbb Q(\sqrt{2})$. Therefore, the degree of the field $\mathbb K$ is equal to $D = 2$.
		Let
		\begin{equation*}
			\Gamma_4 = 1- \dfrac{(d_1-d_210^{m-n})\cdot10^n\cdot\alpha^{-k}}{18}.
		\end{equation*}
		If $\Gamma_4=0$, then $\alpha^{k}\in \mathbb Q$, which is false for $k>0$. By using the properties of the absolute $logarithmic$ $height$, we get
		\begin{equation*}
			h(\gamma_1)=h(\alpha)=\dfrac{\log\alpha}{2},\hspace{0.2cm}
			h(\gamma_2)= \log 10, 
		\end{equation*}
		and 
		\begin{equation*}
			\begin{split}
				h(\gamma_3) &=  h\biggr(\dfrac{(d_1-d_210^{m-n})}{18}\biggr)\\ &\leq h(9)+ h(2)+h(d_1)+h(d_2)+(n-m)\log10 +\log 2 \\ & <7.98+(n-m) \log 10.
			\end{split}
		\end{equation*}
		So we can take $A_1=\log \alpha$, $A_2= 2\log 10$, $A_3= 15.96+2(n-m)\log 10$. Since $n < k+3$, we take $D$ = max$\{n, k, 1\} = k+3$. 
		Thus, taking into account \eqref{19} and using Theorem \ref{thm1}, we obtain
		\begin{equation*}
			3\cdot\alpha^k>|\Gamma_4|> e^{(C\cdot(1+\log 2)(1+\log(k+3))\cdot\log \alpha\cdot2\log 10\cdot(15.96+2(n-m)\log 10))},
		\end{equation*}
		where $C= -1.4\cdot30^6\cdot3^{4.5}\cdot 2^2$. By a simple computation, it follows that 
		\begin{equation}\label{20}
			k\log\alpha-\log 3 <7.9\cdot10^{12}\cdot(1+\log(k+3))(15.96+2(n-m)\log 10).
		\end{equation}
		Using \eqref{16} and \eqref{20}, a computer search with $Mathematica$ gives us that $k < 2.4\cdot10^{28}$.\\

		Now, let us reduce the upper bound on $n$ by using the Baker–Davenport algorithm given in Lemma \ref{lem1}. Let
		\begin{equation*}
			z_3=k\log\alpha - n\log10 +\log\biggr(\dfrac{9}{2d_1}\biggr).
		\end{equation*}
		
		From \eqref{15}, we have
		\begin{equation*}
			|x| = |e^{z_3}-1|< \dfrac{9.81}{10^n-m}<\frac{1}{10}
		\end{equation*}
		for $n-m \geq 2$. Choosing $a = 0.1$, we get the inequality
		\begin{equation*}
			|z_3|=|\log(x+1)|< \dfrac{\log(10/9)}{1/10} \cdot \frac{9.81}{10^{n-m}}<(10.34)\cdot 10^{m-n}
		\end{equation*}
		by Lemma \ref{lem2}. Thus, it follows that 
		\begin{equation*}
			0<\biggr|k\log\alpha - n\log10 +\log\biggr(\dfrac{9}{2d_1}\biggr)\biggr|< (10.34)\cdot 10^{m-n}.
		\end{equation*}
		Dividing this inequality by $\log 10$, we obtain
		\begin{equation}\label{21}
			0<\biggr|k\biggr(\frac{\log\alpha}{\log 10}\biggr) - n +\dfrac{\log\biggr(\dfrac{9}{2d_1}\biggr)}{\log 10}\biggr|< (4.5) \cdot 10^{m-n}.
		\end{equation}
		We can take $\tau=\frac{\log\alpha}{\log 10} \notin \mathbb Q$ and $M = 2.4\cdot 10^{28}$.\\ Then we found that $q_{59} = 808643106803003389273254071835$, the denominator of the $59^{th}$ convergent of $\tau$ is greater than $6M$.
		Now take $\mu=\dfrac{\log(9/2d_1)}{\log 10}$.
		In this case, considering the fact that $1 \leq d_1 \leq9$, a quick computation with $Mathematica$ gives us that $\epsilon(\mu) := \left\lVert \mu q_{59}\right\rVert - M\left\lVert \tau q_{59}\right\rVert = 0.401882$.
		Let $A=4.5$, $B= 10$, and $\omega=n-m$ in Lemma \ref{lem1}. Thus, using $Mathematica$, we can say that \eqref{21} has no solution if 
		\begin{equation*}
			\dfrac{\log(Aq_{59}/\epsilon(\mu))}{\log B} < 31.95<n-m.
		\end{equation*}
		So $n-m \leq 31 $.
		Substituting this upper bound for $n - m$ in \eqref{20}, we obtain $k < 3.1 \cdot 10^{15} $. Now, let
		\begin{equation*}
			z_4= n\log10-k\log\alpha +\log\biggr(\dfrac{d_1-d_210^{m-n}}{18}\biggr).
		\end{equation*}
		From \eqref{19}, we have 
		\begin{equation*}
			|x| = |e^{z_4}-1|< \dfrac{3}{\alpha^k}<\frac{1}{10}
		\end{equation*}
		for $k\geq 25$. Choosing $a = 0.1$, we get the inequality
		\begin{equation*}
			|\Lambda_2|=|\log(x+1)|< \dfrac{\log(10/9)}{1/10} \cdot \frac{3}{\alpha^k}<(3.17)\cdot \alpha^{-k}
		\end{equation*}
		by Lemma \ref{lem2}. Thus, it follows that
		\begin{equation*}
			0<\biggr|n\log10-k\log\alpha +\log\biggr(\dfrac{d_1-d_210^{m-n}}{18}\biggr)|< (3.17)\cdot \alpha^{-k}.
		\end{equation*}
		Dividing both sides by $\log \alpha$, we obtain
		\begin{equation}\label{22}
			0<\biggr|n\biggr(\frac{\log 10}{\log \alpha}\biggr) - k +\dfrac{\log\biggr(\dfrac{d_1-d_210^{m-n}}{18}\biggr)}{\log \alpha}\biggr|< (6.59)\cdot \alpha^{-k}.
		\end{equation}
		Putting $\tau = \frac{\log 10}{\log \alpha} \notin \mathbb Q$ and taking $M= 3.1\cdot10^{15}$, we found that $q_{36}=73257846218558279$, the denominator of the $36^{th}$ convergent of $\tau$ exceeds $6M$. Now take
		$\mu = \dfrac{\log\biggr(\dfrac{d_1-d_210^{m-n}}{18}\biggr)}{\log \alpha}$.
		
		In this case, considering the fact that $1 \leq d_1, d_2 \leq 9$ and $2 \leq n - m \leq 31$, a quick computation gives us the inequality $\epsilon(\mu) := \left\lVert \mu q_{36}\right\rVert - M\left\lVert \tau q_{36}\right\rVert = 0.421589$.
		Let $A=6.59$, $B= \alpha$, and $\omega=k$ in Lemma \ref{lem1}. Thus, using $Mathematica$, we can say that \eqref{22} has no solution if 
		\begin{equation*}
			\dfrac{\log(Aq_{36}/\epsilon(\mu))}{\log B} < 23.58< k.
		\end{equation*}
		This implies $k \leq 23$, which contradicts our assumption that $k\geq 25$. Thus, the proof is completed.

	\end{proof}	
\end{theorem}
{\bf Data Availability Statements:} Data sharing not applicable to this article as no datasets were generated or analyzed during the current study.

{\bf Funding:} The authors declare that no funds or grants were received during the preparation of this manuscript.

{\bf Declarations:}

{\bf Conflict of interest:} On behalf of all authors, the corresponding author states that there is no Conflict of interest.

\end{document}